\documentclass[12pt,a4paper]{article}

\newenvironment{theorem}[1]{\endgraf  {\bf #1}\em}{ \endgraf  }
\def\proof{{\bf Proof. }}
\usepackage{amssymb, amsmath}
\usepackage{graphics}

\DeclareMathOperator{\E}{\sf E}

\DeclareMathOperator{\PP}{\sf P}

\DeclareMathOperator{\Var}{Var}

\begin{document}

\title{\bf\large\MakeUppercase{
New estimates of the convergence rate in the Lyapunov theorem}}

\author{Ilya Tyurin\thanks{Moscow State University, Department of Probability Theory, Moscow. E-mail: itiurin@gmail.com}}
%\date{}

\maketitle

\begin{abstract}
We investigate the convergence rate in the Lyapunov theorem when
the third absolute moments exist. By means of convex analysis we
obtain the sharp estimate for the distance in the mean metric
between a probability distribution and its zero bias
transformation. This bound allows to derive new estimates of the
convergence rate in terms of Kolmogorov's metric as well as the
metrics $\zeta_r$ $(r=1,2,3)$ introduced by Zolotarev. The
estimate for $\zeta_3$ is optimal. Moreover, we show that the
constant in the classical Berry-Esseen theorem can be taken as
$0{.}4785$. In addition, the non-i.i.d. analogue of this theorem
with the constant $0{.}5606$ is provided.
\end{abstract}

Our results \cite{T1} concerning the convergence rate in the
Lyapunov central limit theorem were published in "Doklady Akademii
Nauk" (the article was presented by Professor Yu.~V.~Prokhorov on
June 10, 2009). The complete proofs \cite{T2} were submitted to
the "Theory of Probability and its Applications" on June 8, 2009.
As it turned out later, independently of us Professor Goldstein
has obtained some results that coincide with ours. Namely, an
estimate for the proximity in the mean metric between a
probability distribution and its zero bias transformation, and the
upper bound of the constant in the mean central limit theorem have
been established. His article \cite{G} appeared on arXiv more than
two weeks later, i.e. on June 28, 2009.

The present paper includes not only the results of \cite{T1,T2},
but also their improvements. We show that the constant $C$ in the
Berry-Esseen inequality does not exceed $0{.}4785$. Moreover, we
find a bound for the constant that appears in the generalization
of this theorem in the case of nonidentically distributed
summands. For this case we obtain the estimate $C \leqslant
0{.}5606$. These new results \cite{T3} were presented by Professor
A. V. Bulinski to the "Russian Mathematical Surveys" on November
17, 2009.
\section{Introduction} Consider centered independent (real-valued)
random variables (r.v.) $X_1, \ldots, X_n$ with variances
 $\sigma_1^2,\ldots,\sigma_n^2$ and finite absolute moments $\beta_1,\ldots,\beta_n$.
 We denote
$$
\sigma^2 = \sigma^2 (n) :=\sum_{j=1}^n
\sigma_j^2,\quad\varepsilon_n := \frac{1}{\sigma^3} \sum_{j=1}^n
\beta_j.
$$
According to the Lyapunov theorem,
$S_n:=(X_1+\ldots+X_n)/\sigma(n)$ converges in distribution to the
standard normal r.v. when $\varepsilon_n \to 0$. From both
theoretical and practical points of view, it is very important to
estimate the convergence rate in this theorem. It is known
\cite{Berry,Esseen1} that there exists a minimal numerical
constant $C$ such that for the Kolmogorov distance between $S_n$
and the standard normal variable $N$ holds the inequality
\begin{equation}\label{uniform}
\rho(S_n,N):= \sup\limits_{x\in \mathbb R} \bigl| \PP(S_n
\leqslant x) - \PP(N\leqslant x)\bigr| \leqslant
C\varepsilon_n,\enskip n \in \mathbb N.
\end{equation}
There are plenty of works devoted to estimation of this constant.
Esseen \cite{Esseen1} showed that $C\leqslant 7{.}5$. Bergstr\"om
\cite{Bergstrom} obtained the bound $C\leqslant 4{.}8$. Takano
\cite{Takano} established that in the case of independent
identically distributed (i.i.d.) summands $C \leqslant 2{.}031$.
Zolotarev \cite{Zolot_ner,Zolot1,Zolot2,Zolot_cf} obtained a new
inequality allowing to estimate the proximity of two sums of
independent~r.v. With the help of this inequality he showed
successively that $C \leqslant 1{.}322$ and $C \leqslant
0{.}9051$, while in the case of i.i.d. variables $C \leqslant
1{.}301$ and $C \leqslant 0{.}8197$. The proposed method was
further developed in the works of van Beek \cite{VanBeek} and
Shiganov \cite{Shiganov}, who proved the estimates $C \leqslant
0{.}7975$ and $C \leqslant 0{.}7915$, respectively. For the sums
of identically distributed r.v. Shiganov obtained the bound $C
\leqslant 0{.}7655$, which was sharpened in 2006 by
Shevtsova~\cite{Shevtsova}. She showed that in this case $C
\leqslant 0{.}7056$. In \cite{T1,T2} we derived the estimates $C
\leqslant 0{.}6379$ in the general case and $C \leqslant 0.5894$
for  identically distributed summands. In the present paper we
improve them.

From a private communication with Korolev and Shevtsova we know
that recently they have established the convergence rate in the
central limit theorem in a variety of sences
\cite{KSL,KSBE,KSREF}. In these works only the i.i.d. case was
considered and the bound for the constant $C$ is not as sharp as
ours. However, interesting estimates of the other kind were
obtained.

It is worth mentioning the related problem of determining the
asymptotically best constants in Lyapunov's theorem. As it was
shown by Esseen \cite{Esseen2}, if all the r.v. $X_j,\enskip j =
1, 2, \ldots$ have the same distribution, then
\begin{equation}\label{easympt}
\limsup\limits_{n\to\infty}\frac{\rho(S_n,N)}{\varepsilon_n}
\leqslant C_1:=\frac{\sqrt{10}+3}{6\sqrt{2\pi}}=0{.}409{\ldots},
\end{equation}
and the constant on the right-hand side of this inequality cannot
be lowered (hence the lower bound $C \geqslant C_1$). This result
was elaborated by Rogozin \cite{Rogozin}, who established that
under the same assumptions
\begin{equation}\label{rasympt}
\limsup\limits_{n\to\infty}\frac{\rho(S_n,\mathcal
N)}{\varepsilon_n} \leqslant C_2 := \frac{1}{\sqrt{2\pi}},
\end{equation}
where $\rho(S_n,\mathcal N) := \inf\limits_{G\in\mathcal
N}\rho(S_n,G)$ and $\mathcal N$ is the set of all normal r.v.

Chistyakov \cite{Chi1,Chi2,Chi3} generalized (\ref{easympt}) and
(\ref{rasympt}) to the case of nonidentically distributed
summands. He proved that
\begin{equation*}
\rho(S_n,N) \leqslant C_1 \varepsilon_n +
r_1(\varepsilon_n),\qquad \rho(S_n,\mathcal N) \leqslant C_2
\varepsilon_n + r_2(\varepsilon_n),
\end{equation*}
where $r_1(\varepsilon_n), r_2(\varepsilon_n)$ are
$o(\varepsilon_n)$ when $\varepsilon_n \to 0$.

There are also the estimates of the convergence rate in the
Lyapunov theorem provided that the moments of the order $2+\delta$
exist (see \cite{KS1, KS2}).

Analogues of (\ref{uniform}) are known for other probability
metrics as well, for example, $\zeta_r$ (where $r=1,2,3$). The
latter will be described in detail in section 2. Estimates in
terms of these metrics can be obtained in a natural way using
Stein's method. The proof of the estimates mentioned above uses,
in particular, the so-called zero bias transformation of a
probability distribution (see \cite{zb}).

For the distance in terms of metrics $\zeta_r$ $(r=1,2,3)$ the
following estimates (see \cite{Raic}) are known:
\begin{equation}\label{zetaconst}
\zeta_1(S_n,N) \leqslant 3\varepsilon_n,\qquad \zeta_2(S_n,N)
\leqslant \frac{3\sqrt{2\pi}}{8}\varepsilon_n,\qquad
\zeta_3(S_n,N) \leqslant \frac{1}{2}\varepsilon_n.
\end{equation}

Hoeffding \cite{Hoeffding} considered the problem of finding the
least upper bound of ${\sf E}f(X_1,\ldots,X_n)$ over the set of
all collections of independent simple r.v. satisfying $m$
restrictions of the form $\E g_{ij}(X_j) = c_{ij}, j = 1,\ldots,
n$. More precisely, it was established that in this case one has
to consider only r.v. taking at most $m+1$ values. In the present
work results of \cite{Hoeffding} are generalized to the case of
arbitrary quasiconvex functional defined on the set of all
probability distributions.

The results obtained allowed us to derive an unimprovable estimate
for the proximity in the mean metric between a probability
distribution and its zero bias transformation. The latter was used
to estimate the accuracy of the Gaussian approximation for the
sums of independent variates. It was established that the values
of constants in (\ref{zetaconst}) can be taken 3 times lower. In
addition, our estimate for the metric $\zeta_3$ is optimal.
Furthermore, new estimates for the difference between the
characteristic functions of the normalized sum and the standard
normal r.v. were derived, which allowed us to prove that $C
\leqslant 0{.}5606$ and in the case of i.i.d. summands $C
\leqslant 0{.}4785$.

\section{Main notions and results}  Let $(S, d)$ be a metric
space and denote by $Q$ the set of all finite signed measures on
the Borel $\sigma$-algebra $\mathcal B(S)$ with the operations of
multiplication by a scalar and addition defined as follows: let
$\mu, \mu_1, \mu_2 \in Q, \,\, c \in \mathbb R$, then for each $A
\in \mathcal B(S)$
$$
(c\mu)(A) := c \cdot \mu(A),\,\,\,\,\,\, (\mu_1 + \mu_2)(A) :=
\mu_1(A) + \mu_2(A).
$$
It is easy to see that $Q$ forms a linear space. And the set D of
discrete probability distributions that are concentrated on finite
sets of points is a convex subset of $Q$. The latter means that
$\alpha \mu_1 + (1-\alpha) \mu_2 \in D\,$ for arbitrary $\mu_1,
\mu_2 \in D$ and $\alpha\in(0,1)$.

Consider the set of all collections consisting of $n$ independent
r.v. $X_1, \ldots,X_n$. Then
$$
{\sf E}f(X_1,\ldots,X_n) = \int\limits_{\mathbb R^n} f
dP_{X_1}\ldots dP_{X_n},
$$
where $f:\mathbb R^n \to \mathbb R$, $P_{X_1},\ldots, P_{X_n}$ are
the distributions of $X_1, \ldots, X_n$. Thus, ${\sf
E}f(X_1,\ldots,X_n)$ can be regarded as a function on the set of
measures, which is linear with respect to each of its $n$
arguments.

A function $g: G\to \mathbb R$, where $G$ is a convex set, is said
to be quasiconvex, if for any $x, y \in G$ and $\alpha \in (0,1)$,
we have
$$
g(\alpha x+(1-\alpha) y) \leqslant \max \{g(x), g(y)\}.
$$

We assume that on $S$ some real-valued functions $h_1, \ldots,
h_m$ are defined. Consider the set
$$K := \{ \mu \in D: \,\langle h_i ,\mu \rangle = 0,\,\, i=
1,\ldots,m\},\enskip\text{where}\enskip \langle f, \mu\rangle:=
\int\limits_S fd\mu.$$ It is easy to see that $K$ is convex. Let
$K_j$ be the set of measures $\mu \in K$ that are concentrated on
at most $j$ points ($j \in \mathbb N$).

\begin{theorem}{Theorem 1.} For any quasiconvex function $g: K \to
\mathbb R$, we have
\begin{equation*}
\sup\limits_{\mu \in K} g(\mu) = \sup\limits_{\mu \in K_{m+1}}
g(\mu).
\end{equation*}
\end{theorem}
In this expression, we assume that the supremum over the empty set
is zero.

\begin{theorem}{Theorem 2.} Let $f$ be a nonnegative function on $S$, $V$ -- a linear space with the norm $\|\cdot\|$, $A:K\to V$ --
such a mapping that \begin{equation}\label{lin} A(\alpha
\mu+(1-\alpha)\nu) = \alpha A\mu + (1-\alpha)A\nu
\end{equation}
for arbitrary $\mu, \nu \in K, \alpha \in (0,1)$. Then the least
value of $\gamma$ such that the inequality
\begin{equation} \label{mincnst} \|A\mu\|\leqslant \gamma \langle
f,\mu\rangle
\end{equation}
holds for every measure $\mu \in K$, coincides with the least
value of $\gamma$ such that $(\ref{mincnst})$ is true for every
measure $\mu \in K_{m+1}$.
\end{theorem}

Let $W$ be a zero-mean r.v. with variance $\sigma^2>0$. A r.v.
$W^*$ is said to have the $W$-zero biased distribution if
\begin{equation}\label{zero_bias}
{\sf E} Wf(W) =  \sigma^2 {\sf E}f'(W^*)
\end{equation}
for every differentiable function $f:\mathbb R \to \mathbb R$ such
that the left-hand side of (\ref{zero_bias}) is defined. It is
known (see \cite{zb}) that $W^*$ exists for every $W$ as described
above and has a density
\begin{equation}\label{density}p(w)=\begin{cases}\sigma^{-2}{\sf E}
\left(W\cdot {\bf 1} \{W>w\}\right), & \text{if \,$w\geqslant 0$;}\\
\sigma^{-2}{\sf E} \left(-W\cdot {\bf 1} \{W<w\}\right), &
\text{if \,$w<0$.}
\end{cases}\end{equation}

For every function $f\in C^{(r-1)}(\mathbb R)$, where $r\in\mathbb
N$, define
$$
M_r (f) := \sup\limits_{x \ne y} \left| \frac{f^{(r-1)}(x) -
f^{(r-1)}(y)}{x-y}\right|.
$$
As usual, $C^{(0)}(\mathbb R) := C(\mathbb R)$. If $f\notin
C^{(r-1)}(\mathbb R)$, we set $M_r(f) = \infty$. Denote
\begin{equation*}
\zeta_r (X,Y) := \sup \{ \left| {\sf E}f(X) - {\sf E}f(Y)\right|:
f \in \mathcal F_r \},\enskip r= 1,2,\ldots,
\end{equation*}
where $\mathcal F_r$ is the set of all real bounded functions with
$M_r(f)\leqslant 1$.

Note that $\zeta_1$ has alternative representations. These are the
so-called mean metric
$$
\varkappa_1(X,Y) := \int_{-\infty}^\infty \bigl| \PP(X\leqslant x)
- \PP(Y \leqslant x)\bigr|dx,
$$
and the minimal $L_1$-metric
$$
l_1(X,Y) := \inf\left\{\E|\widetilde X-\widetilde Y|:
Law(\widetilde X) =Law(X), Law(\widetilde Y) = Law(Y)\right\}.
$$
For details see \cite[p. 21]{Zolotarev}.

\begin{theorem}{ Theorem 3.} If $W$ is a centered r.v. with unit
variance and finite third absolute moment, then
\begin{equation}\label{wzb}\zeta_1(W,W^*) \leqslant \frac{1}{2}{\sf E}|W|^3,\end{equation}
with equality when $W$ has a 2-point distribution.
\end{theorem}

\begin{theorem}{Corollary 1.} Consider a r.v. $S_n^*$ having the $S_n$-zero biased distribution.
Then \begin{equation*}\zeta_1(S_n,S_n^*) \leqslant
\frac{1}{2}\varepsilon_n.
\end{equation*}
\end{theorem}
\begin{theorem}{Theorem 4.} The following inequalities are true:
\begin{gather}
\zeta_1(S_n,N) \leqslant 2\zeta_1(S_n,S_n^*) \leqslant
\varepsilon_n,\qquad \zeta_2(S_n,N) \leqslant
\frac{\sqrt{2\pi}}{4}
\zeta_1(S_n,S_n^*)\leqslant \frac{\sqrt{2\pi}}{8}\varepsilon_n,\\
\label{zetaimproved} \zeta_3(S_n,N) \leqslant \frac{1}{3}
\zeta_1(S_n,S_n^*) \leqslant \frac{1}{6}\varepsilon_n.
\end{gather}
The latter double inequality is optimal, namely, for every $\delta
> 0$ there exists such a sequence of i.i.d. r.v. $X_1, X_2,
\ldots$, that
\begin{equation*}\frac{\zeta_3(S_n,N)}{\varepsilon_n} \geqslant
\frac{1}{6}-\delta,\quad n=1,2,\ldots\end{equation*}
\end{theorem}
For $\gamma > 0$ and $t \in \mathbb R$ we set
$$
b(t,\gamma):= \begin{cases}  - t^2 + 2\gamma a |t|^3, & \text{if $\gamma |t| < M$;}\\
-2 \left(\frac{1}{\gamma}\right)^2\left(1 - \cos{\gamma t}
\right), & \text{if $M
\leqslant \gamma |t|\leqslant 2 \pi$;}\\
0, & \text{if $\gamma |t| > 2 \pi$.}
\end{cases}
$$
Here $$a:= \max\limits_{x>0} \{( \cos(x) - 1 +
{x^2}/{2})/{x^3}\}\approx 0{.}099162,$$ and $M$ is the point where
this maximum is attained, $M \approx 3{.}995896$.

Denote $f_{S_n}(t) := {\sf E}e^{itS_n}$, $\varphi(t) :=
\exp\left(-t^2/2\right)$, $\delta_n(t):=|f_{S_n}(t)-\varphi(t)|$,
$t \in \mathbb R$.
\begin{theorem}{Theorem 5.} For every
$t\in \mathbb R$ we have
\begin{gather}\label{n10} |f_{S_n}(t)| \leqslant \widehat f_1(\varepsilon_n,t) := \exp\left(\frac{1}{2}b( t,
2\varepsilon_n)\right), \\
\label{n105}\delta_n(t) \leqslant
\widehat\delta_1(\varepsilon_n,t) :=
\varepsilon_n\varphi(t)\int\limits_0^{|t|}
\frac{s^2}{2}\exp\left(\frac{s^2}{2}\right) ds.
\end{gather}
Define $A:=\varepsilon_n^{-1/3}/6a$. For all $t\in\mathbb R$ the
following estimate is true
\begin{equation}\label{n12}
\delta_n(t) \leqslant\widehat\delta_2(\varepsilon_n,t)
:=\begin{cases} \varepsilon_n \varphi(t)\int\limits_0^{|t|} \frac{
s^2}{2}\exp\left(\frac{s^2\varepsilon_n^{2/3}}{2}\right) ds, &
\text{$|t| \leqslant A$;}\\\varepsilon_n\varphi(t)
\left(\int\limits_0^A
\frac{s^2}{2}\exp\left(\frac{s^2\varepsilon_n^{2/3}}{2}\right) ds
+ \int\limits_A^{|t|}\frac{s^2}{2 l}\exp\left( 2 a \varepsilon_n
s^3\right)ds \right), & \text{$|t| > A$.}
\end{cases}
\end{equation}
where \begin{equation*}l := \inf\limits_{t \geqslant
0}\left\{\exp\left(-{t^2/2}+2a t^3)\right)\right\}\approx
0{.}624489.\end{equation*}
\end{theorem}
The quantities $\widehat\delta_1(\varepsilon, t)$ and
$\widehat\delta_2(\varepsilon, t)$ can be expressed in terms of
the so-called Dawson integral
$$
{\rm Daw}(t):=\exp\left(-t^2\right)\int\limits_0^t
\exp\left(s^2\right)ds,
$$
which can be computed by the means of several efficient numerical
procedures. For example, such a function is available in the GNU
Scientific Library (GSL). It is easy to check that
$$
\widehat\delta_1(\varepsilon, t) = \varepsilon\left(\frac{t}{2} -
\frac{1}{\sqrt{2}} {\rm Daw}\left(\frac{t}{\sqrt{2}}\right)
\right).
$$
Moreover,
\begin{equation*}
\widehat\delta_2(\varepsilon, t) = \begin{cases}
\exp\left(\frac{t^2(\varepsilon^{2/3}-1)}{2}\right)
\left(\frac{t\varepsilon^{1/3}}{2} - \frac{1}{\sqrt{2}}{\rm
Daw}\left(\frac{t\varepsilon^{1/3}}{\sqrt{2}}\right)\right), &
\text{$|t|
\leqslant A$;}\\
\exp\left(\frac{(1/6a)^2-t^2}{2}\right)\left(\frac{(1/6a)^2}{2} -
\frac{1}{\sqrt{2}}{\rm Daw}(1/6a\sqrt{2}) \right) +
\frac{\varphi(t)}{12al}\exp\left(2a\varepsilon
u^3\right)\mid_{u=A}^t\, & \text{$|t|
> A$.}
\end{cases}
\end{equation*}
These representations are of great importance, since they allow to
reduce significantly the amount of numerical calculations required
for the proof of Theorem 7.

In the case of i.i.d. variables the estimates can be slightly
improved. Denote $\tau_n~:=~\frac{1}{\sigma^3}\sum_{j=1}^n
\sigma_j^3$, and let $X_1, X_2, \ldots$ be centered i.i.d. r.v.
with unit variances and finite third absolute moments $\beta$.
Then $\varepsilon_n = \beta/\sqrt{n}$ and $\tau_n = 1/\sqrt{n}$.
\begin{theorem}{Theorem 6.} For the sequence of r.v. as defined above and every $t\in \mathbb R$
\begin{gather}\label{n7} |f_{S_n}(t)| \leqslant \widehat f_2(\varepsilon_n, n, t):=\left(1 + \frac{b( t,
\varepsilon_n+1/\sqrt{n})}{n}\right)^{n/2},\\
\label{n8}\delta_n(t) \leqslant
\widehat\delta_3(\varepsilon_n,n,t):=
\varepsilon_n\varphi(t)\int\limits_0^{|t|} {\left(1 + \frac{b( s,
\varepsilon_n+1/\sqrt{n})}{n}\right)^\frac{n-1}{2} \frac{s^2}{2}
\exp\left(\frac{s^2}{2}\right)}\,ds.
\end{gather}
Let $m\in \mathbb N$ and $n \geqslant m$. Then
\begin{gather}
\label{n85} |f_{S_n}(t)| \leqslant \widehat
f_3(\varepsilon_n,m,t):=\exp\left(\frac{1}{2}b( t,
\varepsilon_n + 1/\sqrt{m})\right),\\
\label{n9}\delta_n(t) \leqslant
\widehat\delta_4(\varepsilon_n,m,t):=
\varepsilon_n\varphi(t)\int\limits_0^{|t|} \exp\left(\frac{m-1}{2
m}\,b( s, \varepsilon_n+1/\sqrt{m}) + \frac{s^2}{2}\right)
\frac{s^2}{2}\,ds.
\end{gather}
\end{theorem}
Estimates (\ref{n10})-(\ref{n9}) allowed to establish the
following result.
\begin{theorem}{Theorem 7.} The constant $C$ in inequality $(\ref{uniform})$ does not exceed $0{.}5606$, and in the case of
identically distributed summands $C\leqslant 0{.}4785$.
\end{theorem}
\section{Proofs}{\bf Proof of Theorem 1.} If $K =
\varnothing$, then $K_{m+1} = \varnothing$, and the statement of
our theorem is true. Further we suppose that the set $K$ is
nonempty.

The sequence of the sets $K_1, K_2,\ldots$ increases to the set
$K$. Therefore,
$$
\sup\limits_{\mu \in K} g(\mu) = \sup\limits_{j\geqslant 1}
\sup\limits_{\mu \in K_j} g(\mu) = \sup\limits_{j\geqslant m+1}
\sup\limits_{\mu \in K_j} g(\mu).
$$
So, it is sufficient to show that
$$\sup\limits_{\mu \in K_{m+1}} g(\mu)\geqslant\sup\limits_{\mu
\in K_{m+2}} g(\mu)\geqslant\sup\limits_{\mu \in K_{m+3}}
g(\mu)\geqslant\ldots$$ Let's take an arbitrary measure $\mu\in
K_j$, where $j>m+1$, and show that there exists $\mu'\in K_{j-1}$
such that $g(\mu')\geqslant g(\mu)$.

Let $\mu$ be concentrated in points $s_1, \ldots, s_j \in S$ and
\begin{equation}\label{sootvetstvie}\mu(s_i) = \mu_i
\geqslant 0, \enskip i = 1, \ldots, j.\end{equation} The vector
$\bar{\mu} = (\mu_1, \ldots, \mu_j)$ defines a probability
distribution, so
\begin{equation}\label{norm}\mu_1 + \ldots + \mu_j = 1.\end{equation}
Moreover, the conditions  $\,\langle h_i,\mu\rangle = 0,\enskip i=
1,\ldots,m,\,$ hold and therefore
\begin{equation}\label{matrix}
\mu_1 \cdot h_i(s_1)+ \ldots + \mu_j \cdot h_i(s_j) = 0,\quad i=
1,\ldots,m.
\end{equation}

Vice versa, an arbitrary vector with nonnegative coordinates
satisfying the system of linear equations (\ref{norm}) and
(\ref{matrix}) defines according to (\ref{sootvetstvie}) an
element of the set $K_j$, and if one of its coordinates equals
zero -- an element of $K_{j-1}$. We have $m+1$ equations and at
least $m+2$ unknowns, so there exists a nonzero solution
$\bar{\nu} = (\nu_1,\ldots,\nu_j)$ of the corresponding
homogeneous system. Since the sum of the coordinates of this
vector is equal to zero, but the vector itself is nonzero, it
follows that $\bar{\nu}$ has both positive and negative
coordinates. Therefore, there exist the least $\alpha\geqslant 0$
and the least $\beta \geqslant 0$ such that one of the coordinates
of the vector $\bar{\mu_*} = \bar{\mu} - \alpha \bar{\nu}$ equals
zero and some coordinate of $\bar{\mu^*} = \bar{\mu} + \beta
\bar{\nu}$ is equal to zero. If $\alpha = 0$, then $\mu\in
K_{j-1}$. Otherwise,
$$\bar{\mu} = \frac{\beta}{\alpha + \beta}\bar{\mu_*} + \frac{\alpha}{\alpha +
\beta}\bar{\mu^*},$$ and because of the quasiconvexity
$g(\mu)\leqslant \max\{ g(\mu_*),g(\mu^*)\}$, where $\mu_*, \mu^*$
are the distributions defined by $\bar{\mu_*}$ and $\bar{\mu^*}$.
Thus, $g(\mu)\leqslant g(\mu_*)$ or $g(\mu)\leqslant g(\mu^*)$.
But $\mu_*$ and $\mu^* \in K_{j-1}.\enskip\square$

{\bf Proof of Theorem 2.} According to Theorem 1, it is sufficient
to prove that for every fixed value of $\gamma$ the function
$$g(\mu) := \|A\mu\|- \gamma \langle f,\mu\rangle$$ is quasiconvex. Let $\alpha+\beta = 1$.
By the properties of the norm
\begin{multline*}
\|A(\alpha \mu+\beta \nu)\| -  \gamma \langle f,\alpha \mu+\beta
\nu\rangle = \| \alpha A\mu + \beta A\nu\|-  \gamma \langle
f,\alpha \mu\rangle -  \gamma \langle f,\beta \nu\rangle \leqslant
\\ \leqslant \| \alpha A\mu\| + \|\beta A\nu\|-  \gamma \langle
f,\alpha \mu\rangle -  \gamma \langle f,\beta \nu\rangle =\alpha
g(\mu)+ \beta g(\nu) \leqslant \max\{g(\mu),g(\nu)\}.
\enskip\square
\end{multline*}

{\bf Proof of Theorem 3.} We begin by showing that without loss of
generality we can consider simple r.v. $W$. It is sufficient to
establish that for every r.v. $W$ satisfying the conditions of the
theorem there exists a sequence $(W_n)_{n\geqslant 1}$ of simple
r.v. with zero means and unit variances such that
\begin{equation}\label{need}\zeta_1(W_n,W_n^*)\to\zeta_1(W,W^*) \enskip\text{and}\enskip
{{\sf E} |W_n|^3}\to {{\sf E} |W|^3},\quad n\to
\infty.\end{equation} We suppose that r.v. $W$ is defined on the
probability space $(\mathbb R, \mathcal B(\mathbb R), \PP = P_W)$
and construct a sequence of simple r.v.
$\left(W'_n\right)_{n\geqslant 1}$ that converges to the r.v. $W$
in $L_3$ norm. We set
$$W_n:= \frac{ W'_n - {\sf E}W'_n}{\sqrt{\Var W'_n}}.$$ It is easy to see
that $W_n$ converges to $W$ in $L_3$ as well. Therefore, the
second condition in (\ref{need}) is obviously satisfied. It
remains to show that the first one also holds.

From the triangle inequality for the metric $\zeta_1$ one can
easily derive that
\begin{equation}\label{difference}
|\zeta_1(W,W^*) - \zeta_1(W_n,W_n^*)|\leqslant \zeta_1(W,W_n) +
\zeta_1(W^*,W_n^*).
\end{equation}
The first summand on the right-hand side of (\ref{difference})
tends to zero, since
$$
\zeta_1(W,W_n) = l_1(W,W_n)\leqslant \E |W-W_n|\leqslant \left(\E
|W-W_n|^3\right)^{\frac{1}{3}}.
$$
Let's evaluate the second summand. For the function $f\in\mathcal
F_1$ we set $F(x) := \int_0^x f(u) du$. Then
\begin{equation}\label{newdif}{\sf E}f(W^*) - {\sf E}f(W_n^*) =
{\sf E} W F(W) - {\sf E} W_n F(W_n).\end{equation} The difference
of expectations on the left-hand side of (\ref{newdif}) does not
change if we replace the function $f(x)$ by $f(x) - f(0)$.
Therefore, we can assume without loss of generality that $f(0) =
0$. Then $|f(x)-f(y)|\leqslant |x-y|$ yields $|f(x)|\leqslant |x|$
and so $|F(x)|\leqslant |x|^2, |x f(x)|\leqslant |x|^2$. According
to the finite-increment theorem,
$$W F(W) - W_n F(W_n) = (W-W_n) \cdot \{x
F(x)\}'\mid_{x=\xi}=(W-W_n) \{F(\xi) + \xi f(\xi)\},$$ where $\xi$
is a number between $W$ and $W_n$. Moreover,
$$|F(\xi)+\xi f(\xi)| \leqslant
2|\xi|^2 \leqslant 2 (|W| + |W_n|)^2.$$ This gives the estimate
$$|{\sf E}\{W F(W) - W_n F(W_n)\}|\leqslant 2 {\sf E}|W - W_n| \left( |W| + |W_n| \right)^2.$$
And finally, the H\"older's inequality yields
$${\sf E}|W - W_n| \left( |W| + |W_n|\right)^2 \leqslant \left({\sf E}|W-W_n|^3\right)^{\frac{1}{3}} \left({\sf E}\left(|W| + |W_n|\right)^{3}\right)^{\frac{2}{3}}.$$
Obviously, $\left({\sf E}|W-W_n|^3\right)^{\frac{1}{3}}\to 0$,
since $W_n$ converges to $W$ in $L_3$. Thus, the second summand in
(\ref{difference}) tends to zero and (\ref{need}) is fulfilled.

So, it is sufficient to consider simple r.v. Let $A_1$ be a
function that  maps the distribution $P_X$ of a r.v. $X$ to its
zero-biased distribution $P_X^*$. Moreover, consider a linear
operator $A_2$ that maps a signed measure $\nu$ to its cumulative
distribution function (c.d.f.)
${G_\nu(x):=\nu\left((-\infty,x]\right)}$. It is easy to see that
$$
(\alpha P_{W_1} + (1-\alpha) P_{W_2})^* = \alpha P_{W_1}^* +
(1-\alpha) P_{W_2}^*,
$$ hence the mapping $A_2 - A_2 A_1$ satisfies
(\ref{lin}). If we set $h_1(x) = x$, $h_2(x) = x^2-1$ and apply
Theorem 2 to $f(x) = |x|^3$, $A = A_2 - A_2 A_1$, $V$ -- the
normed space of integrable functions on the real line with the
norm
$$
\| G \|= \int_{-\infty}^\infty |G(x)| dx,
$$
then the problem reduces to the case of simple r.v. taking at most
3 values. C.d.f. of a simple r.v. $W$ is a staircase function.
Using formula (\ref{density}) one can easily obtain the c.d.f. of
$W^*$. Therefore, it is not difficult to find the explicit
expression for $\varkappa_1(W,W^*)$.
\begin{center}{\includegraphics{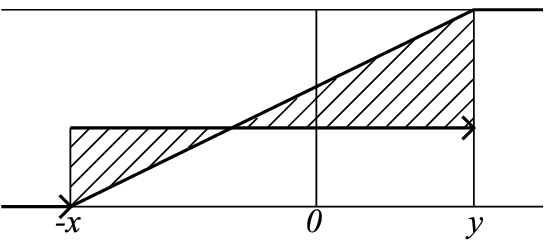} \qquad
\raisebox{0.06cm}{\includegraphics{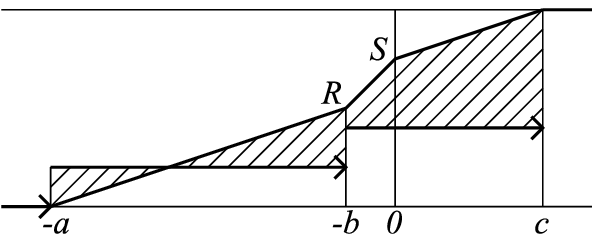}}}

{Pic. 1}
\end{center}

Let $W$ take exactly two values $-x$ and $y$ with probabilities
$p$ and $q$. Then its c.d.f. is piecewise constant and has two
steps in points $-x$ and $y$ that are equal to $p$ and $q$,
respectively. Since $W$ is centered, we have $px = qy$, which
together with (\ref{density}) yields that $W^*$ is uniformly
distributed on $[-x,y]$. Therefore, on $[-x,y]$ its c.d.f. is
linear, and its graph is a segment that connects $(-x,0)$ and
$(y,1)$. By definition $\varkappa_1(W,W^*)$ equals the area of the
figure bounded by distribution functions of these r.v. (in the
case considered it is a union of two triangles, see pic. 1, left).

It follows from the conditions $\E W =0$ and $\E W^2= 1$ that
$x=\sqrt{q/p},\enskip y = \sqrt{p/q}.$ Hence
\begin{equation}\label{beta3}
\E |W|^3 = p x^3 + q y^3 = q \sqrt{\frac{q}{p}} + p
\sqrt{\frac{p}{q}}.
\end{equation}
Let's find the area of the figure bounded by c.d.f. of r.v. $W$
and $W^*$. Density of $W^*$ is $px = qy = \sqrt{pq}$. Thus, the
slope of the c.d.f. of this r.v. on [-x,y] is $\sqrt{pq}$. The
length of the vertical leg of the first triangle is $p$, and that
of the second one is $q$. Hence, the total area of both triangles
is
$$
\frac{1}{2} p^2 \frac{1}{\sqrt{pq}} + \frac{1}{2} q^2
\frac{1}{\sqrt{pq}} = \frac{1}{2} p \sqrt{\frac{p}{q}}  +
\frac{1}{2} q \sqrt{\frac{q}{p}} = \frac{1}{2} \E |W|^3.
$$
Therefore, if $W$ takes exactly two values, there is equality in
(\ref{wzb}).

Consider the case when $W$ takes three values. We assume without
loss of generality that two of them ($-a < -b$) do not exceed zero
and one ($c$) is positive. As before, the c.d.f. of $W$ is
piecewise constant and the c.d.f. of $W^*$ is piecewise linear.
However, the form of the figure bounded by them is more
complicated
 (see pic. 1, right). Denote by $R$ the value of the c.d.f. of
$W^*$ at $-b$ and by $S$ its value at $0$. Let $W$ take values
$-a, -b, c$ with probabilities $p, q, r$, respectively. Then,
because of the moment-type restrictions,
$$
\left\{
  \begin{array}{l}
    p + q + r = 1, \\
    -p a -q b + rc = 0, \\
    p a^2 +q b^2 + r c^2 =1.\\
  \end{array}
\right.
$$
It is a system of linear equations with respect to $p,q,r$. Using
Cramer's rule, we obtain
$$
p = (1-bc)(c+b)/\Delta, \qquad q = (ac-1)(a+c)/\Delta, \qquad r =
(ab+1)(a-b)/\Delta,
$$
where $\Delta = (a+c)(b+c)(a-b)$. Thus, every r.v. with zero mean
and variance 1 that takes three values is uniquely determined by
these three values. It is easy to see that $p,q,r$ are nonnegative
iff
\begin{equation}\label{acbc}ac \geqslant 1, \qquad bc \leqslant 1.\end{equation}
In other words, a r.v. $W$ taking the values $-a, -b, c$ exists
iff (\ref{acbc}) is satisfied. Our aim is to prove that the
function
\begin{equation}\label{defg}g(a,b,c):=\varkappa_1(W,W^*) -
\frac{1}{2}\E|W|^3\end{equation} does not exceed zero. Its
explicit form in terms of variables $a,b,c$ depends on how the
c.d.f. of r.v. $W$ and $W^*$ are located with respect to each
other. There are 5 cases:

\noindent I. $R \leqslant p$, $S \leqslant p + q$, or,
equivalently, $a(a-b) \leqslant 1$, $c \geqslant 1$. In this case
$$
g(a,b,c) = \frac{r}{c} - p a^3 - q b^3.
$$
II. $R \leqslant p$, $S \geqslant p + q \Leftrightarrow a(a-b)
\leqslant 1$, $c \leqslant 1$.
$$
g(a,b,c) = \frac{r}{c} - p a^3 - q b^3.
$$
III. $R \geqslant p$, $S \leqslant p + q \Leftrightarrow a(a-b)
\geqslant 1$, $c \geqslant 1$ (the latter implies $c(b+c)\geqslant
1$).
$$
g(a,b,c) = p a \left\{a-b-\frac{1}{a} \right\}^2 + \frac{r}{c} -
pa^3 -qb^3.
$$
IV. $p\leqslant R \leqslant p + q$, $S \geqslant p + q
\Leftrightarrow a(a-b) \geqslant 1$, $c(b+c) \geqslant 1$, $c
\leqslant 1$.
$$
g(a,b,c) = p a \left\{a-b-\frac{1}{a} \right\}^2 + \frac{r}{c} -
pa^3 -qb^3.
$$
V. $R \geqslant p+q \Leftrightarrow c(b+c) \leqslant 1$.
$$
g(a,b,c) = \frac{p}{a} - r c^3.
$$
Note that in each of these cases $g$ is the same function defined
by (\ref{defg}). As a result, if the values $a,b,c$ satisfy the
restrictions of two cases simultaneously, then for the function
$g$ we can use the expression corresponding to any of them.

As one can see, in the cases I, II, and in the cases III, IV the
function $g$ has the same representation. Therefore, further we
distinguish three possibilities:

\noindent A. $a(a-b) \leqslant 1$.
$$
g(a,b,c) = \frac{r}{c} - p a^3 - q b^3.
$$
\noindent B. $a(a-b) \geqslant 1$, $c(b+c) \geqslant 1$
$$
g(a,b,c) = p a \left\{a-b-\frac{1}{a} \right\}^2 + \frac{r}{c} -
pa^3 -qb^3.
$$
\noindent C. $c(b+c) \leqslant 1$.
$$
g(a,b,c) = \frac{p}{a} - r c^3.
$$

We show that in each of the cases A, B and C the function $g$ does
not exceed zero.

\noindent Case A. \begin{multline*}g(a,b,c) = \frac{r}{c} - p a^3
- q b^3 = \frac{(1+ab)(a-b)}{\Delta c} -
\frac{(1-bc)(b+c)a^3}{\Delta} - \frac{(ac-1)(a+c)b^3}{\Delta} = \\
= \frac{(a-b)(ac-1)(1-bc)}{\Delta c}(-1-bc-ac-ab) \leqslant
0\end{multline*} because of (\ref{acbc}), nonnegativeness of
$a,b,c$ and the fact that $a > b$.

\noindent Case C. \begin{multline*}g(a,b,c) = \frac{p}{a} - r c^3
=
\frac{(1-bc)(b+c)}{\Delta a} - \frac{(1+ab)(a-b)c^3}{\Delta} = \\
= \frac{(ac-1)}{\Delta a} \left\{-\left[(ab+1)(a-b)\right]c^2 -
\left[1+ b(a-b)\right] c - b\right\}.\end{multline*}

Since $ac \geqslant 1$, it suffices to prove that the expression
enclosed by braces does not exceed zero. Consider this expression
as a function of the variable $c$ while holding the others fixed.
When $c = 0$, this function equals $-b\leqslant 0$. Moreover, it
decreases with respect to $c$, since the coefficients of terms $c$
and $c^2$ are negative. Consequently, for all positive values of
$c$ it does not exceed zero.

\noindent Case B. \begin{multline*}g(a,b,c) = p a
\left\{a-b-\frac{1}{a} \right\}^2 + \frac{r}{c} - pa^3 -qb^3 =
\frac{(1-bc)(b+c)a}{\Delta} \left\{a-b-\frac{1}{a} \right\}^2 + \\
+ \frac{(1+ab)(a-b)}{\Delta c} - \frac{(1-bc)(b+c)a^3}{\Delta} -
\frac{(ac-1)(a+c)b^3}{\Delta} = \frac{1-bc}{\Delta ac}\left\{ k_2
c^2 + k_1 c + k_0\right\},
\end{multline*}
where $\enskip k_2 = 1-2a(a-b)(1+ab),\enskip k_1 = b\left[ 1 -
a(a-b)(1+ab)\right],\enskip k_0 = a(a-b)(1+ab).$ Assume that
$g(a,b,c)
> 0$. Due to the condition $a(a-b) \geqslant 1$ one has
\begin{gather*} 1- a(a-b)(1+ab) \leqslant 1-
(1+ab) = -ab \leqslant 0, \\
1 - 2a(a-b)(1+ab) \leqslant 1 - 2(1+ab) = -1 -2ab < 0.
\end{gather*}
Therefore, $k_2$ and $k_1$ do not exceed zero and, consequently,
$k_2 c^2 + k_1 c + k_0$ decreases with respect to $c$.
Consequently, if one reduces the value of the variable $c$ while
holding $a$ and $b$ fixed, $g$ will remain positive. The variable
$c$ is bounded from below by two conditions:
$$ac \geqslant 1\enskip \text{and} \enskip c(b+c)\geqslant 1.$$
The first of these conditions can be omitted, since it follows
from the other two:
$$c(b+c)\geqslant 1\enskip \text{and} \enskip a(a-b)\geqslant 1.$$
Indeed, let $ac < 1$. Then
$$
1 \leqslant c(b+c) < \frac{1}{a}\left( b + \frac{1}{a} \right)
\Rightarrow a^2 < ab + 1 \Rightarrow a(a-b) < 1.
$$
Therefore, we can reduce $c$ to the value $c_*$ such that
$c_*(b+c_*) = 1$. And $g$ will remain positive. But the situation,
when $c(b+c) = 1$, satisfies the restrictions of the case C, for
which we established that $g \leqslant 0$ -- a contradiction.
$\enskip\square$

{\bf Proof of Corollary 1.} Without loss of generality assume
$\sigma$ = 1. Let $I$ be a random index taking values $1,\ldots,
n$ with probabilities $\sigma_1^2,\ldots,\sigma_n^2$, independent
of $X_1,\ldots,X_n$. Construct on an extended probability space
$$
S_i' := \sum_{j\ne i} X_j + X_i^*, \quad i = 1,\ldots, n,
$$
where $X_i^*$ has the $X_i$-zero biased distribution and is
independent of $I, X_1,\ldots,X_n$, ${i=1,\ldots,n}$. Then
$S_n^*=S_I'$ (see~\cite{zb}). Therefore, for an arbitrary function
$f\in \mathcal F_1$ one has
\begin{equation}\label{snzvezda} {\sf E}f(S_n^*) ={\sf E}f( S_I') = \sum\limits_{k=1}^n
{\sf E} f(S_I') {\bf 1}\{I = k\}= \sum\limits_{k=1}^n {\sf E}
f(S_k') {\bf 1}\{I = k\} = \sum\limits_{k=1}^n \sigma_k^2 {\sf E}
f(S_k').
\end{equation}
Consequently,
\begin{multline*}
|{\sf E}f(S_n) - {\sf E}f(S_n^*)| = \left|\sum\limits_{k=1}^n
\sigma_k^2 {\sf E} f(S_n) - \sum\limits_{k=1}^n \sigma_k^2{\sf
E}f(S_k')\right|\leqslant \sum\limits_{k=1}^n \sigma_k^2
\left|{\sf E} f(S_n) - {\sf E}f(S_k')\right|\leqslant\\
\leqslant\sum\limits_{k=1}^n \sigma_k^2 \zeta_1(X_k, X_k^*) =
\sum\limits_{k=1}^n \sigma_k^3 \zeta_1\left(\frac{X_k}{\sigma_k},
\frac{X_k^*}{\sigma_k}\right) \leqslant\sum\limits_{k=1}^n
\sigma_k^3 \cdot\frac{1}{2}\E\left|\frac{X_k}{\sigma_k}\right|^3 =
\frac{1}{2}\sum\limits_{k=1}^n \beta_k = \frac{1}{2}\varepsilon_n.
\end{multline*}
Here we used the statement of Theorem 3 for r.v.
$\frac{1}{\sigma_k}X_k$ as well as the property of homogeneity of
the metric $\zeta_1$ (i. e. $\zeta_1(c X,c Y) = c \zeta_1(X,Y)$)
and $(\alpha X)^* \stackrel{D}=\alpha X^* $. $\enskip\square$

{\bf Proof of Theorem 4.} It is easy to see that for continuous
$f: \mathbb R \to \mathbb R$
$$
h(w) := e^{\frac{w^2}{2}}\int\limits_{-\infty}^w(f(w)-{\sf
E}f(N))e^{-\frac{x^2}{2}}dx
$$
satisfies the Stein's equation
$$
h'(w) - w h(w) = f(w) - {\sf E}f(N).
$$
Hence \begin{equation*} |{\sf E}f(S_n) - {\sf E}f(N)| = |{\sf
E}h'(S_n) - {\sf E}S_n h(S_n)| = |{\sf E} h'(S_n) - {\sf
E}h'(S_n^*)| \leqslant M_2(h) \zeta_1(S_n,S_n^*).
\end{equation*}
As it was shown in \cite{Raic},
$$
M_2(h) \leqslant \min\{ 2 M_1(f), \frac{\sqrt{2 \pi}}{4}M_2(f),
\frac{1}{3}M_3(f)\}.
$$
Taking into account that $M_r(f)\leqslant 1$ in the definition of
the metric $\zeta_r$, one has
\begin{equation}\label{zetabiased}\zeta_1(S_n,N) \leqslant
2\zeta_1(S_n,S_n^*),\quad \zeta_2(S_n,N) \leqslant
\frac{\sqrt{2\pi}}{4} \zeta_1(S_n,S_n^*),\quad \zeta_3(S_n,N)
\leqslant \frac{1}{3} \zeta_1(S_n,S_n^*).\end{equation} Estimates
in terms of $\varepsilon_n$ are obtained by applying Corollary 1
to (\ref{zetabiased}).

Let's prove the optimality of (\ref{zetaimproved}). We set $f(x)
:= x^3/6$. Then $M_3(f) = 1$ and ${\E f(N) = 0}$, since the r.v.
$N$ is symmetric and the function $f$ is odd. Consider $X_1,
X_2,\ldots$ -- a sequence of i.i.d. variables with zero means and
unit variances. Then
$$
\E S_n ^3 = \frac{\E X_1^3}{\sqrt{n}}  \quad\text{and}\quad
\varepsilon_n = \frac{\E |X_1|^3}{\sqrt{n}}.
$$
As a result, we have $$\frac{\zeta_3(S_n, N)}{\varepsilon_n}
\geqslant |\E f(S_n)- \E f(N)| = \frac{1}{6} \frac{|\E X_1^3|}{\E
|X_1|^3}.$$ It only remains to prove that $|\E X_1^3|/\E |X_1|^3$
can be arbitrarily close to unity. According to (\ref{beta3}), the
third absolute moment of a centered r.v. $X_1$ with variance 1
taking two values $x=-\sqrt{q/p},\,\,y = \sqrt{p/q}$ with
probabilities $p$ and $q$, respectively, is equal to $$ \E|X_1|^3
= q\sqrt{\frac{q}{p}} + p \sqrt{\frac{p}{q}}.
$$
It's easy to see that the third moment of this r.v. equals
$$ \E X_1^3 = - q \sqrt{\frac{q}{p}} + p \sqrt{\frac{p}{q}}.
$$
Obviously, $|\E X_1|/\E|X_1|^3 \to 1$ when $p\to 1$. $\enskip
\square$

\begin{theorem}{Lemma 1 {\rm (\cite{weitere})}.} Let $W$ be centered r.v. with variance 1
and finite third absolute moment~$\beta$. Denote $f(t):= {\sf
E}e^{itW}, t \in \mathbb R$. Then for all $t \in \mathbb R$
\begin{equation}\label{pervoe} |f(t)|^2 \leqslant 1 +
b(t,\beta+1),
\end{equation}
moreover,
\begin{equation}
\label{noniid} |f_{S_n}(t)|^2 \leqslant \exp\left(
b(t,\varepsilon_n+\tau_n)\right).
\end{equation}
\end{theorem}
\begin{theorem}{Lemma 2.} For every $t \in \mathbb R$ the function $b(t,\gamma)$ is nondecreasing with respect to
$\gamma$.
\end{theorem}
\proof This can be checked directly by calculating the
derivative.$\enskip\square$
\begin{theorem}{Lemma 3.} If $W$ is a centered r.v. with variance 1 and $f(t) = \E e^{itW}$, then
\begin{equation}\label{intestimate}
|f(t) - \varphi(t)| \leqslant \varphi(t)\int\limits_{0}^t |f(s) -
f^*(s)| s \exp\left( \frac{s^2}{2}\right) ds, \enskip t \in
\mathbb R,
\end{equation}
where $f^*(t)$ is the characteristic function of a r.v. $W^*$
having the $W$-zero biased distribution.
\end{theorem}
\proof According to the definition of the $W$-zero biased
distribution,
\begin{equation}\label{complexbias}
f'(t) = i\E W\cos(tW) -\E W\sin(tW) = -i t\E \sin(tW^*)
-t\E\cos(tW^*)  = -t \E e^{itW^*}.
\end{equation}
Consider the function $\psi(t):=\frac{f(t)}{\varphi(t)}$. Note
that $\psi(0) = 1$. Taking into account (\ref{complexbias}), we
have
$$
\psi'(t) = \frac{d}{dt}\left(f(t) e^{\frac{1}{2}t^2}\right) =
f'(t) e^{\frac{1}{2}t^2} + t f(t)  e^{\frac{1}{2}t^2} = \{f(t) -
f^*(t)\}t e^{\frac{1}{2}t^2}.
$$
Then
$$
\left| \frac{f(t)}{\varphi(t)} - 1\right| = |\psi(t) - \psi(0)|
\leqslant \int_{0}^t |\psi'(s)| ds = \int_0^t|f(s) - f^*(s)| s
e^{\frac{1}{2}s^2}ds. \enskip\square
$$
\begin{theorem}{Lemma 4.}
For arbitrary r.v. $X$ and $Y$ we have
$$
\left|\E e^{itX} - \E e^{itY}\right| \leqslant t\zeta_1(X,Y).
$$
\end{theorem}
\proof It is well known that for all $t, x, y\in \mathbb R$ holds
the inequality
$$
\left|e^{itx} - e^{ity}\right| \leqslant t|x-y|.
$$
Thus, for arbitrary $\widetilde X, \widetilde Y$ defined on one
probability space such that $Law(\widetilde X) = Law(X)$ and
$Law(\widetilde Y) = Law(Y)$, we have
\begin{equation}\label{long}\left|\E e^{itX} - \E e^{itY}\right|=\left|\E e^{it\widetilde
X} - \E e^{it\widetilde Y}\right| \leqslant\E\left|
e^{it\widetilde X} -  e^{it\widetilde Y}\right|  \leqslant
t\E|\widetilde X-\widetilde Y|.\end{equation} Passing in
(\ref{long}) to the greatest lower bound among every possible
$\widetilde X, \widetilde Y$, we obtain
$$
\left|\E e^{itX} - \E e^{itY}\right| \leqslant t l_1(X,Y) =
t\zeta_1(X,Y). \enskip\square
$$

{\bf Proof of Theorem 5.} The inequality (\ref{n10}) is a
consequence of Lemma 1. Indeed, according to the Lyapunov
inequality, we have $\sigma_j^3 \leqslant \beta_j, j = 1,\ldots,
n$. Hence $\tau_n \leqslant \varepsilon_n$. Now (\ref{n10})
follows from (\ref{noniid}) and Lemma 2.

Further we assume without loss of generality that $\sigma = 1$.
Denote $f_j(t) := \E e^{itX_j},\,f_j^*(t) := \E
e^{itX_j^*},\,\,j=1,\ldots,n,$ and set $W:=S_n$ in
(\ref{intestimate}). Using Lemma 4 and Corollary~1 we get
\begin{equation*}
\left|f(s) - f^*\left(s\right)\right| \leqslant
s\zeta_1(S_n,S_n^*)\leqslant \frac{\varepsilon_n s}{2}.
\end{equation*}
Substituting the latter into (\ref{intestimate}), we arrive at
(\ref{n105}).

According to (\ref{snzvezda}),
\begin{equation*}
f^*\left(s\right) = \sum\limits_{m=1}^n \sigma_m^2 {\sf E} e^{is
S_m'} =\sum\limits_{m=1}^n \sigma_m^2 f_m^*(s)\prod\limits_{j\ne
m} f_j(s).
\end{equation*}
Therefore,
\begin{multline}\label{fsfs}
\left|f(s) - f^*\left(s\right)\right| = \left| \sum\limits_{i=1}^n
\sigma_i^2 \prod\limits_{j} f_j(s) - \sum\limits_{i=1}^n
\sigma_i^2 f_i^*(s)\prod\limits_{j\ne i} f_j(s)\right| = \\ =
\left|\sum\limits_{i=1}^n \sigma_i^2 \left\{ f_i(s) -
f_i^*(s)\right\}\prod\limits_{j\ne i} f_j(s)\right|.
\end{multline}
From Lemma 4 and Theorem 3 we have \begin{equation}\label{fjs}
|f_j(s) - f_j^*(s)| \leqslant s \zeta_1(X_j,X_j^*) = s \sigma_j
\zeta_1\left(\frac{X_j}{\sigma_j},\frac{X_j^*}{\sigma_j}\right)
\leqslant \frac{\beta_j s}{2\sigma_j^2}.
\end{equation}
It follows from (\ref{pervoe}) that $|f_j(s)| \leqslant
\exp(-{\sigma_j^2 s^2/2}+2\beta_j a |s|^3)$ for all real $s$. As a
result,
\begin{multline}\label{first}\left|f(s) - f^*\left(s\right)\right|\leqslant
\left|\sum\limits_{i=1}^n \beta_i \frac{s}{2} \prod\limits_{j\ne
i} \exp\left(-\frac{\sigma_j^2 s^2}{2}+2\beta_j a s^3\right)\right| = \\
= \left|\sum\limits_{j=1}^n \beta_j \frac{s}{2}
\exp\left(\frac{\sigma_j^2 s^2}{2}-2\beta_j a s^3\right)\right|
\exp\left(-\frac{ s^2}{2}+2\varepsilon_n a
s^3\right).\end{multline} Since $\sigma_j^3 \leqslant \beta_j
\leqslant \varepsilon_n$ for $j=1,\ldots,n$, we have for such $j$
\begin{multline}\label{second}
\exp\left(\frac{\sigma_j^2 t^2}{2}-2\beta_j a t^3\right) \leqslant
\exp\left(\frac{\sigma_j^2 t^2}{2}-2\sigma_j^3 a t^3\right) =
\left.\exp\left(\frac{s^2}{2}-2 a
s^3\right)\right|_{s=t\sigma_j}\leqslant \\ \leqslant \sup
\left\{\exp\left(\frac{s^2}{2}-2 a s^3\right):
s\in[0,t\varepsilon_n^{1/3}] \right\}.
\end{multline}
The function $\exp\left(\frac{s^2}{2}-2 a s^3\right)$ increases on
the segment $[0, 1/6a]$ and at the point $1/6a$ it attains its
global maximum equal to $1/l$. Therefore,
\begin{equation}\label{third}\sup \left\{\exp\left(\frac{s^2}{2}-2
a s^3\right):
s\in[0,t\varepsilon_n^{1/3}] \right\} = \begin{cases}\exp\left(\frac{t^2 \varepsilon_n^{2/3}}{2}-2 a\varepsilon_n t^3\right), & \text{if \,$t\varepsilon_n^{1/3} \leqslant 1/6a$;}\\
1/l, & \text{otherwise.}
\end{cases}\end{equation}
Combining (\ref{first}), (\ref{second}) and (\ref{third}) gives
for $s\varepsilon_n^{1/3} \leqslant 1/6a$
$$
\left|f(s) - f^*\left(s\right)\right|\leqslant \frac{\varepsilon_n
s}{2}\exp\left(\frac{s^2\varepsilon_n^{2/3}}{2} -\frac{s^2}{2}
\right), $$ and for $s\varepsilon_n^{1/3} > 1/6a$
$$
\left|f(s) - f^*\left(s\right)\right|\leqslant\frac{\varepsilon_n
s}{2 l}\exp\left( -\frac{s^2}{2} + 2 a \varepsilon_n s^3 \right).
$$
Substituting the expressions obtained into (\ref{intestimate}), we
get the required estimates.

{\bf Proof of Theorem 6.} At first we prove (\ref{n7}). Denote
$f_1(t):=\E e^{itX_1}$. According to Lemma 1,
\begin{equation}\label{bylo}
\left|f_1\left(\frac{t}{\sqrt{n}}\right)\right|\leqslant \sqrt{1 +
b\left(\frac{t}{\sqrt{n}}, \beta + 1 \right)} = \sqrt{1 +
\frac{1}{n}b\left(t, \frac{\beta+1}{\sqrt{n}} \right)} = \sqrt{1 +
\frac{1}{n}b\left(t, \varepsilon_n + \tau_n \right)}.
\end{equation}
Now (\ref{n7}) follows from the fact that
$f_{S_n}(t)=f_1^n(t/\sqrt{n})$.

To establish (\ref{n85}) we note that $1+ x \leqslant e^x$ for all
real $x$. Applying this inequality to (\ref{n7}) gives
\begin{equation}\label{form} |f_{S_n}(t)| \leqslant
\exp\left(\frac{1}{2} b( t, \varepsilon_n + \tau_n)\right).
\end{equation}
It remains to note that the sequence $(\tau_m)_{m\geqslant 1}$
decreases, which leads to (\ref{n85}).

We set $W:=S_n$ in Lemma 3. Applying (\ref{fsfs}) to the r.v.
$\frac{1}{\sqrt{n}}X_1, \ldots, \frac{1}{\sqrt{n}}X_n$ yields
$$|f(s) - f^*(s)| = \left|f_1\left(\frac{s}{\sqrt{n}}
\right)\right|^{n-1}\cdot\left| f_1\left(\frac{s}{\sqrt{n}}
\right) - f_1^*\left(\frac{s}{\sqrt{n}} \right)\right|.$$ The
first factor can be estimated with the help of (\ref{bylo}) and
the second -- by means of (\ref{fjs}). We have
\begin{equation}\label{lastestimate} |f(s) - f^*(s)| \leqslant
\left(1 + \frac{b( s,
\varepsilon_n+\tau_n)}{n}\right)^\frac{n-1}{2} \frac{\varepsilon_n
s}{2}.
\end{equation}
Substituting the expression obtained into (\ref{intestimate}), we
get (\ref{n8}). To establish (\ref{n9}) we apply the inequality
${1+x\leqslant e^x}$ to the first factor on the right-hand side of
(\ref{lastestimate}) and note that $(\tau_m)_{m\geqslant 1}$ is
decreasing. $\enskip\square$

{\bf Proof of Theorem 7.} Let $D(\varepsilon,n)$ denote the least
quantity such that for every collection consisting of $n$ r.v.
$X_1,\ldots, X_n$ with $\varepsilon_n = \varepsilon$ holds the
inequality $$\rho(S_n,N) \leqslant
D(\varepsilon,n)\varepsilon_n.$$ We set $$D(\varepsilon) :=
\sup\limits_{n\geqslant 1} {D(\varepsilon,n)}.$$ Then the constant
$C$ can be determined as $$C = \sup_{\varepsilon
> 0} D(\varepsilon).$$ Hence, it suffices to show that for all possible values of $\varepsilon$ and
$n$ the quantity ${D(\varepsilon,n)\leqslant 0{.}5606}$ (and in
the case of i.i.d. r.v. $D(\varepsilon,n)\leqslant 0{.}4785$). For
$\varepsilon\geqslant 1/0{.}5606$ (respectively,
${\varepsilon\geqslant 1/0{.}4785}$) the latter is obvious, since
$\rho(S_n,N)\leqslant 1$.

Moreover, denote
$\lambda_n:=\sigma^2(n)/(\sigma^2(n)-\max\limits_{k=1,\ldots,n}\sigma_k^2)$
and set
$$\widehat\varepsilon_n:=\lambda_n^{3/2}\varepsilon_n,\quad
\varepsilon'_n:=\lambda_n^{3/2}\tau_n,\quad
\varepsilon''_n:=\lambda_n^2\sum_{k=1}^n\sigma^4_k/\sigma^4(n).$$
Then, according to the inequality (I.52) from \cite{PraCLT}, for
$\widehat\varepsilon_n+\varepsilon'_n\leqslant 0{.}2$ we have
\begin{equation}\label{I52}\rho(S_n,N)\leqslant 0{.}27283\widehat\varepsilon_n +
0{.}19948\varepsilon'_n+0{.}09116\varepsilon''_n+0{.}00095(\widehat\varepsilon_n+\varepsilon'_n)^2.\end{equation}

Assume without loss of generality that $\sigma = 1$. Then the
Lyapunov inequality yields $\sigma_j^3 \leqslant \beta_k \leqslant
\sum_{k=1}^n\beta_k= \varepsilon_n, \,j= 1,\ldots,n$. Hence
$\lambda_n\leqslant(1-\varepsilon_n^{2/3})^{-1}$. In addition,
$\varepsilon'_n\leqslant\widehat\varepsilon_n$ and, as it was
shown in \cite{PraCLT},
$\varepsilon''_n\leqslant(\varepsilon'_n)^{4/3}$. From these
inequalities and (\ref{I52}) it follows easily that
$D(\varepsilon)\leqslant 0{.}5606$ when $\varepsilon\leqslant
0{.}02$.

In the case of i.i.d. summands $\varepsilon_n =
\beta_1/(\sigma_1^3\sqrt{n}) \geqslant 1/\sqrt{n}$. Thus,
$n\geqslant \lceil 1/\varepsilon_n^2\rceil$ and
\begin{equation}\label{lam}\lambda_n=\frac{n}{n-1}\leqslant
\frac{n_0(\varepsilon_n)}{n_0(\varepsilon_n)-1},\end{equation}
where $n_0(\varepsilon):=\lceil 1/\varepsilon^2\rceil.$ Moreover,
\begin{equation}\label{varep}\varepsilon'_n = \frac{\lambda_n^{3/2}}{\sqrt{n}},
\quad \text{and} \quad \varepsilon''_n =
\frac{\lambda^2}{n}.\end{equation} Combining (\ref{I52}),
(\ref{lam}) and (\ref{varep}) yields $D(\varepsilon)\leqslant
0{.}4785$ when $\varepsilon\leqslant 0{.}037$. Therefore, in the
general case we have to consider $\varepsilon$ from the segment
$I_1 =\left[0{.}02;{1}/{0{.}5606}\right]$ and in the case of
i.i.d. r.v. -- from $I_2=\left[0{.}037;{1}/{0{.}4785}\right]$. The
proof for these values of $\varepsilon$ is based on an inequality
due to Prawitz \cite{distr}
\begin{multline}\label{praineq}
\frac{\rho(S_n,N)}{\varepsilon_n}\leqslant\frac{1}{\varepsilon_n}\left(\int\limits_{-U_0}^{U_0}\frac{1}{U}\left|K\left(\frac{u}{U}\right)\right|\cdot|\delta_n(u)|du+
\int\limits_{U_0<|u|\leqslant
U}\frac{1}{U}\left|K\left(\frac{u}{U}\right)\right|\cdot|f_n(u)|du+\right.
\\
\left.+\int\limits_{-U_0}^{U_0}\left|\frac{1}{U}K\left(\frac{u}{U}\right)-\frac{i}{2\pi
u}\right|\cdot|\varphi(u)|du+\int\limits_{|u|>U_0}\left|
\frac{\varphi(u)}{2\pi u}\right|du\right),
\end{multline}
where $\enskip K(u) := \frac{1}{2}(1-|u|) +
\frac{i}{2}\left((1-|u|)\cot(\pi
u)+\frac{\textrm{sgn}(u)}{\pi}\right),\quad 0<U_0\leqslant U$.

It follows from (\ref{praineq}) that $D(\varepsilon)$ does not
exceed the quantity $D^*(\varepsilon, U_0, U)$, which arises on
the right-hand side of (\ref{praineq}) when we substitute
$\delta_n(t)$ with its estimate $\min\{\widehat
\delta_1(\varepsilon,t),\widehat \delta_2(\varepsilon,t)\}$,
$|f_n(t)|$ -- with the estimate $\widehat f_1(\varepsilon,t)$ and
select such parameters $U_0, U$ that the resulting expression was
as little as possible. This procedure was carried out with the aid
of computer for several hundreds values of $\varepsilon$ dispersed
on the segment $I_1$. To obtain the estimates for the intermediate
points we used the following property of the quantities $D^*$,
which holds due to the monotonicity of the functions $\widehat
f_1,\ldots,\widehat f_3$ and $\widehat \delta_1,\ldots,\widehat
\delta_4$ with respect to their first arguments.
\begin{equation}\label{deps} D^*(\varepsilon^{(1)},U_0, U)\leqslant
\frac{\varepsilon^{(2)}}{\varepsilon^{(1)}}D^*(\varepsilon^{(2)},U_0,
U), \enskip \varepsilon^{(1)}<\varepsilon^{(2)}.\end{equation}

The extremal value of the quantity $D^*(\varepsilon, U_0, U) =
0{.}56054$ was attained for ${\varepsilon = 0{.}5092}, U_0 =
2{.}4852, U = 5{.}9508.$

In the case of i.i.d. r.v. the estimates were constructed in a
different way.

For the fixed value of $\varepsilon$ we estimated the quantities
$D(\varepsilon, n), \, n\geqslant 1,$ separately . For $n<m$,
where $m$ is some natural number, the individual estimates of
$D(\varepsilon, n)$ were given. On the right-hand side of
(\ref{praineq}) we substituted $\delta_n(t)$ and $|f_n(t)|$ with
their upper estimates $\widehat\delta_3(\varepsilon,n,t)$ and
$\widehat f_2(\varepsilon,n,t)$. After that the computational
procedure as described above was carried out to select the optimal
parameters $U$ and $U_0$. For $n\geqslant m$ the quantities
$D(\varepsilon, n)$ were estimated uniformly. On the right-hand
side of (\ref{praineq}) the estimates $\widehat
\delta_4(\varepsilon,m,t)$ and $\widehat f_3(\varepsilon,m,t)$
were used. As before, it was sufficient to carry out the
calculations only for the finite number of points, since a
property similar to (\ref{deps}) holds in this case as well. For
the i.i.d. r.v. the extremal value $0{.}47849$ was attained for
$\varepsilon = 0{.}3536, n = 8$, $U_0 = 2.6157, U = 8.9115.$

Thus, the constant $C$ does not exceed $0{.}5606$. And if we
restrict to the case of i.i.d. r.v., we have $C\leqslant
0{.}4785$.

\section*{Acknowledgement}
\par The author would like to thank Professor A. V. Bulinski for
useful discussions and valuable advice.

\end{document}